\documentclass[10]{article}
\usepackage{amsfonts}
\usepackage{theorem}

\newtheorem{theorem}{Theorem}
\newtheorem{corollary}{Corollary}
\newtheorem{lemma}{Lemma}

\newtheorem{example}{Example}
\usepackage[unicode,colorlinks,plainpages=false]{hyperref}
\newcommand{\openbox}{$\begin{array}{c}
\hspace*{-0.55em}\sqcap \hspace*{-0.60em}\\[-0.4em] \hline
\multicolumn{1}{c}{\hspace*{-0.60em}}\\[-0.8em]
\end{array}
$}

\usepackage{graphics}
\usepackage{graphicx}
\usepackage{tikz}

\usepackage{enumerate}

\begin{document}

\centerline{\bf A Note on Semigroup Algebras of Permutable Semigroups}
%\centerline{\bf the ideals of a semigroup algebra ${\mathbb F}[S]$}

\bigskip

\centerline{A. Nagy and M. Zubor}

\begin{abstract}

Let $S$ be a semigroup and $\mathbb F$ be a field.
For an ideal $J$ of the semigroup algebra ${\mathbb F}[S]$ of $S$ over $\mathbb F$, let $\varrho _J$ denote the restriction
(to $S$) of the congruence on ${\mathbb F}[S]$ defined by the ideal $J$.
 A semigroup $S$ is called a permutable semigroup if $\alpha \circ \beta =\beta \circ \alpha$ is satisfied for all congruences $\alpha$ and $\beta$ of $S$. In this paper we show that if $S$ is a semilattice or a rectangular band then
$\varphi _{\{S;{\mathbb F}\}}:\ J\mapsto \varrho _J$
is a homomorphism of the semigroup $(Con ({\mathbb F}[S]);\circ )$ into the relations semigroup $({\cal B}_S; \circ )$ if and only if $S$ is a permutable semigroup.
\end{abstract}

\medskip

For a non-empty set $X$, let ${\cal B}_X$ denote the set of all binary relation on a set $X$. It is known that
${\cal B}_X$ form a semigroup under the operation $\circ$ defined by the followong way:
for arbitrary $\alpha ,\beta \in {\cal B}_X$ and $x, y\in X$,
$(x,y)\in \alpha \circ \beta$ if and only if there is an element $z\in X$ such that
$(x,z)\in \alpha$ and $(z,y)\in \beta$.
An algebraic structure $A$ will said to be permutable if the congruences of $A$ commute with
each other (under the operation $\circ$). For arbitrary equivalences
$\alpha$ and $\beta$ of a set $X$, $\alpha \circ \beta$ is an equivalence on $X$ if and only if
$\alpha \circ \beta =\beta \circ \alpha$. Thus an algebraic stucture $A$ is permutable
if and only if the set $Con(A)$ of all congruences on $A$ form a subsemigroup in the
relation semigroup $({\cal B}_A;\circ)$.

Let $S$ be a semigroup and $\mathbb F$ be a field.
%As the semigroup algebra ${\mathbb F}[S]$ is permutable,
%the semilattice
%$(Con({\mathbb F}[S]);\circ )$ is isomorphic to the $\vee$-semilattice
%of all ideals of ${\mathbb F}[S]$; in this paper we identify these two semilattices.
For an ideal $J$ of the semigroup algebra ${\mathbb F}[S]$, let $\varrho _J$
denote the congruence on $S$ which is the restriction of the congruence on ${\mathbb F}[S]$ defined by $J$. For an arbitrary
congruence $\alpha$ on $S$, let $\mathbb{F}[\alpha]$ denote the
kernel of the extended canonical homomorphishm
$\mathbb{F}[S]\rightarrow \mathbb{F}[S/\alpha]$. By Lemma 5 of Chapter 4
of \cite{Okni}, for every semigroup $S$ and every field $\mathbb F$, the mapping
$\varphi_{\{S;{\mathbb F}\}}:\ J \mapsto \varrho _J$ is a surjective homomorphism
of the semilattice $(Con({\mathbb F}[S]);\wedge)$ onto the semilattice
$(Con(S);\wedge)$ such that
$\varrho _{{\mathbb F}[\alpha ]}=\alpha$ for every congruence $\alpha$ on $S$. As a homomorphic image of a semigroup is also a semigroup, and
$\alpha \circ \beta =\alpha \vee \beta$ is satisfied for every congruences $\alpha$ and $\beta$ of a permutable semigroup,
%and a semigroup $S$ is permutable if and only if $Con (S)$ is a subsemigroup of the relation semigroup $({\cal B}_S;\circ )$,
the assertions of the following lemma are obvious.

\begin{lemma}\label{xxx} If $S$ is a semigroup such that, for a field
${\mathbb F}$,
$\varphi _{\{S;{\mathbb F}\}}:\ J\ \to \varrho _J$ is a homomorphism of the
semigroup
$(Con({\mathbb F}[S]);\circ)$ into the semigroup $({\cal B}_S;\circ )$
then $S$ is a permutable semigroup. Moreover, if $S$ is a permutable
semigroup then
$\varphi_{\{S;{\mathbb F}\}}$ is a homomorphism of
$(Con({\mathbb F}[S]);\circ )$ onto the semigroup $(Con(S);\circ )$ if and only if
$\varphi_{\{S;{\mathbb F}\}}$ is a homomorphism of the semilattice
$(Con({\mathbb F}[S]);\vee )$ onto the lattice $(Con(S);\vee )$, that is, $ker_{\varphi_{\{S;{\mathbb F}\}}}$ is $\vee$-compatible.
\end{lemma}

\medskip

The next example shows that the converse of the first assertion of Lemma~\ref{xxx} is not true, in general; for a permutable semigroup $S$, the condition "$\varphi _{\{S;{\mathbb F}\}}$ is a homomorphism
of $(Con({\mathbb F}[S]);\circ )$ onto the semigroup $(Con(S);\circ )$" depends on the field $\mathbb{F}$.

\begin{example}
Let $C_4$, ${\mathbb F}_3$ and $\mathbb{F}_2$ denote the cyclic group of order $4$, the $3$-element field, and the $2$-element field,
respectively. It is known that every group is a permutable semigroup.
Denote the elements of $C_4$ by $1,a,a^2,a^3$ (1 is the identity element).
It is easy to see that
$I=\mathrm{Span}(1+a^2,a+a^3)$ and $J=\mathrm{Span}(1+a,a+a^2,a^2+a^3)$
are ideals of $\mathbb{F}_3[C_4]$, and
\[\varphi_ {\{C_4;{\mathbb F}_3\}}(I)\vee\varphi_ {\{C_4;{\mathbb F}_3\}}(J)=\]
\[\varrho _I\vee \varrho _J=\iota\vee \{(1,a^2),(a,a^3)\}\ne \omega_{C_4}=\varrho _{(I+J)}=\]
\[=\varrho _{(I\vee J)}=\varphi_ {\{C_4;{\mathbb F}_3\}}(I\vee J).\]
Thus $ker_{\varphi_{\{C_4;{\mathbb F}_3\}}}$ is not $\vee$-compatible and so $\varphi_{\{C_4;{\mathbb F}_3\}}$ is not a homomorphism of $Con({\mathbb F}_3[C_4]);\circ)$
onto the semigroup $(Con(C_4);\circ)$.

%However $\varphi_{\{C_4;{\mathbb F}_2\}}$ is a homomorphism of $Con({\mathbb F}_2[C_4]);\circ)$
%onto the semigroup $(Con(C_4);\circ)$.
By computing the pricipal ideals of $\mathbb{F}_2[C_4]$ we get that the principal ideals of $\mathbb{F}_2[C_4]$ form a chain. As the dimension of $\mathbb{F}_2[C_4]$ is finite, from this it follows that $\mathbb{F}_2[C_4]$ has no other ideals. So $Con(\mathbb{F}_2[C_4])$ is the next (here $\alpha_{C_2}$ denotes the congruence of $C_4$ defined by $C_2$):
\begin{center}
\begin{tikzpicture}[scale=1,cap=round,
angle/.style={font=\fontsize{7}{7}\color{black}\ttfamily}
]
	
	\coordinate (a) at (0,2);
	\coordinate (b) at (0,1);
	\coordinate (c) at (0,0);
	\coordinate (d) at (0,-1);
	\coordinate (e) at (0,-2);

	\draw (a) -- (e);
			
	\filldraw[fill=white, draw=white]  (a) circle [radius=0.3];
	\filldraw[fill=white, draw=white]  (b) circle [radius=0.3];
	\filldraw[fill=white, draw=white]  (c) circle [radius=0.3];
	\filldraw[fill=white, draw=white]  (d) circle [radius=0.3];
	\filldraw[fill=white, draw=white]  (e) circle [radius=0.3];
		
	\draw (a) node{$\mathbb{F}_2[C_4]$};
	\draw (b) node{$\mathbb{F}_2[\omega_{C_4}]$};	
	\draw (c) node{$\mathbb{F}_2[\alpha_{C_2}]$};
	\draw (d) node{$Span(1+a+a^2+a^3)$};
	\draw (e) node{$\{0\}$};		
\end{tikzpicture}

\end{center}
It is easy to see that $ker_{\varphi _{\{C_4;{\mathbb F}_2\}}}$ is $\vee$-compatible and so
$\varphi_{\{C_4;{\mathbb F}_2\}}$ is a homomorphism of $Con({\mathbb F}_2[C_4]);\circ)$
onto the semigroup $(Con(C_4);\circ)$.
\end{example}

\medskip
By Lemma~\ref{xxx} and the above Example, it is a natural idea to find all couples $(S,{\mathbb F})$ of permutable semigroups $S$ and fields $\mathbb F$, for which the mapping $\varphi _{\{S;{\mathbb F}\}}$ is a homomorphism of the semigroup
$(\mathrm{Con}({\mathbb F}[S]), \circ )$ onto the semigroup $(Con (S); \circ )$.
In the next we show that if $S$ is an arbitrary permutable semilattice or an arbitrary permutable rectangular band then, for an arbitrary field $\mathbb F$, $\varphi _{\{S;{\mathbb F}\}}$ satisfies the previous condition.

For notations and notions not defined here, we refer to \cite{Clifford, Nagy, Okni, Petrich}

\begin{theorem}\label{semilattice} Let $S$ be a permutable semilattice. Then, for an arbitrary field $\mathbb F$,
$\varphi _{\{S;{\mathbb F}\}}$ is a homomorphism of the semigroup
$(\mathrm{Con}({\mathbb F}[S]), \circ )$ onto the semigroup $(Con (S); \circ )$.
%if and only if $S$ is permutable.
\end{theorem}

\noindent
{\bf Proof}.
%By Lemma~\ref{xxx}, it is sufficient to show that
%$ker _{\varphi _{\{S;{\mathbb F}\}}}$ is $\vee$-compatible
%if $S$ is a permutable semilattice and $\mathbb F$ is an arbitrary field.
Assume that $S$ is a permutable semilattice. Then, by Lemma 2 of \cite{Hamilton},
$|S|\leq 2$.
We can consider the case when $|S|=2$.
Let \[S=\{ e, f\}\quad (e\neq f).\] Then
\[e^2=e,\quad \hbox{and}\quad f^2=f.\]
We can suppose that \[ef=fe=e.\]
It is clear that $S$ has two congruences: $\iota _S$ and $\omega _S$.

Let $\mathbb F$ be an arbitrary field.
It is easy to see that \[J_e=\{ \alpha e:\ \alpha \in {\mathbb F}\}\quad \hbox{and}\quad J_{e-f}=\{\alpha (e-f):\ \alpha \in {\mathbb F}\}\]
are proper ideals of
${\mathbb F}[S]$. As
\[dim(J_e)=dim(J_{e-f})=1,\] the ideals $J_e$ and $J_{e-f}$
are minimal ideals of ${\mathbb F}[S]$.
We show that the ideals of ${\mathbb F}[S]$ are
\[\{ 0\}, J_e, J_{e-f}\quad \hbox{and}\quad {\mathbb F}[S].\]

Let $J\neq \{ 0\}$ be a proper ideal of ${\mathbb F}[S]$.
Clearly $dim(J)=1$. Let \[A=\alpha e+\beta f\in J\] be a non-zero element. Then
$$(\alpha+\beta )e=e(\alpha e+\beta f)=eA\in J.$$
If $\alpha+\beta\ne 0$ then $J=J_e$. If $\alpha+\beta=0$ then $J=J_{e-f}$.

Thus the ideals of ${\mathbb F}[S]$ are
$\{ 0\}$, $J_e$, $J_{e-f}$, ${\mathbb F}[S]$. So $Con(\mathbb{F}[S])$ is the next:

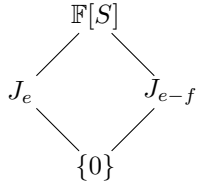
\begin{center}
\begin{tikzpicture}[scale=1,cap=round,
angle/.style={font=\fontsize{7}{7}\color{black}\ttfamily}
]
	
	\coordinate (a) at (0,-1);
	\coordinate (b) at (1,0);
	\coordinate (c) at (0,1);
	\coordinate (d) at (-1,0);
	
	\draw (a) -- (b);
	\draw (a) -- (d);
	\draw (c) -- (b);
	\draw (d) -- (c);
	
	\filldraw[fill=white, draw=white]  (a) circle [radius=0.3];
	\filldraw[fill=white, draw=white]  (b) circle [radius=0.3];
	\filldraw[fill=white, draw=white]  (c) circle [radius=0.3];
	\filldraw[fill=white, draw=white]  (d) circle [radius=0.3];
		
	\draw (a) node{$\{0\}$};
	\draw (b) node{$J_{e-f}$};	
	\draw (c) node{$\mathbb{F}[S]$};
	\draw (d) node{$J_e$};

\end{tikzpicture}
 % \caption{}

\end{center}

It is a matter of checking to see that the $ker_{\varphi_{\{S;{\mathbb F}\}}}$-classes
of $Con({\mathbb F}[S])$ are
$\{ \{ 0\}, J_e\}$ and $\{ J_{e-f}, {\mathbb F}[S]\}$.
It is easy to see that $ker_{\varphi _{\{S;{\mathbb F}\}}}$ is $\vee$-compatible and so, by Lemma~\ref{xxx},
$\varphi _{\{S;{\mathbb F}\}}$ is a homomorphism of the semigroup
$(\mathrm{Con}({\mathbb F}[S]), \circ )$ onto the semigroup $(Con (S); \circ )$.
\hfill\openbox

\begin{corollary}\label{corsemilattice} Let $S$ be a semilattice. Then, for a field $\mathbb F$,
$\varphi _{\{S;{\mathbb F}\}}$ is a homomorphism of the semigroup
$(\mathrm{Con}({\mathbb F}[S]), \circ )$ into the relation semigroup $({\cal B}_S;\circ)$
if and only if $|S|\leq 2$.
\end{corollary}

\noindent
{\bf Proof}. By Lemma~\ref{xxx}, Theorem~\ref{semilattice} and Lemma 2 of \cite{Hamilton}, it is obvious.\hfill\openbox

\begin{theorem}\label{rectangularband} Let $S=L\times R$ be
a permutable rectangular band ($L$ is a left zero semigroup, $R$ is a right zero semigroup).
Then, for an arbitrary field $\mathbb F$, $\varphi _{\{S;{\mathbb F}\}}$ is
a homomorphism of the semigroup
$(\mathrm{Con}({\mathbb F}[S]), \circ )$ onto the semigroup $(Con (S); \circ )$.
%if and only if $S$ is permutable.
\end{theorem}

\noindent
{\bf Proof}.
Let $\mathbb F$ be an arbitrary field and $S=L\times R$ be a permutable rectangular band. As a rectangular band satisfies the identity $axyb=ayxb$, that is, every rectangular band is a medial semigroup, Corollary 1.2 of \cite{Cherubini} implies $|L|\leq 2$ and $|R|\leq 2$.

First consider the case when $|L|=1$. Then $S$ is isomorphic to the right zero
semigroup $R$, and $|S|\leq 2$. We can suppose that $|S|=2$.
Let $S=\{ e, f\}$ ($e\neq f$).
The congruences of $S$ are $\iota _S$ and $\omega _S$. We show that the
ideals of ${\mathbb F}[S]$ are
\[\{ 0\}, J_{e-f}={\mathbb F}[\omega _S]=\{\alpha (e-f):\ \alpha \in {\mathbb F}\}\quad \hbox{and}\quad{\mathbb F}[S].\]

Let $J\neq \{ 0\}$ be an arbitrary ideal.
Assume that there is an element \[0\neq \alpha _0e+\beta _0f\in J\] for which
\[\alpha _0+\beta _0\neq 0\] is satisfied.
Then \[(\alpha _0+\beta _0)e=(\alpha _0e+\beta _0f)e\in J\] and so
\[e={1\over \alpha _0+\beta _0}(\alpha _0+\beta _0)e\in J\] from which we get
$f=ef\in J$. Consequently \[J={\mathbb F}[S].\]

Next, consider the case when \[\alpha _0+\beta _0=0\] is satisfied
for every \[A=\alpha _0e+\beta _0f\in J.\] Then $\beta =-\alpha _0$
and so \[A=\alpha _0e+\beta _0f=\alpha _0e-\alpha _0f=\alpha _0(e-f)\in J_{e-f}.\]
Consequently, \[J\subseteq J_{e-f}.\]
As $dim(J_{e-f})=1$, the ideal $J_{e-f}$ is minimal. Hence \[J=J_{e-f}.\]
Thus the ideals of ${\mathbb F}[S]$ are $\{ 0\}$, $J_{e-f}$ and ${\mathbb F}[S]$, indeed. So $Con(\mathbb{F}[S])$ is

\begin{center}
\begin{tikzpicture}[scale=1,cap=round,
angle/.style={font=\fontsize{7}{7}\color{black}\ttfamily}
]
	
	\coordinate (a) at (0,2);
	\coordinate (b) at (0,1);
	\coordinate (c) at (0,0);

	\draw (a) -- (c);
			
	\filldraw[fill=white, draw=white]  (a) circle [radius=0.3];
	\filldraw[fill=white, draw=white]  (b) circle [radius=0.3];
	\filldraw[fill=white, draw=white]  (c) circle [radius=0.3];
		
	\draw (a) node{$\mathbb{F}[S]$};
	\draw (b) node{$J_{e-f}$};	
	\draw (c) node{$\{0\}$};
		
\end{tikzpicture}
 % \caption{}
\end{center}

It is a matter of checking to see that the $ker_{\varphi _{\{S;{\mathbb F}\}}}$-classes of $Con({\mathbb F}[S])$ are
$\{ \{ 0\} \}$ and $\{ J_{e-f}, {\mathbb F}[S]\}$.
It is easy to see that $ker_{\varphi _{\{S;{\mathbb F}\}}}$ is $\vee$-compatible and so, by Lemma~\ref{xxx},
$\varphi _{\{S;{\mathbb F}\}}$ is a homomorphism of the semigroup
$(\mathrm{Con}({\mathbb F}[S]), \circ )$ onto the semigroup $(Con (S); \circ )$.

If $|R|=1$ then $S$ is a left zero semigroup, and $|S|\leq 2$. We can prove, as in the previous part of the proof,
that $ker_{\varphi _{\{S;{\mathbb F}\}}}$ is $\vee$-compatible.

Next, consider the case when $|L|=|R|=2$. Let\[L=\{ a_1, a_2\},\quad R=\{ b_1, b_2\}.\] Let $\alpha _L$
and $\alpha _R$ denote the kernels of the projection homomorphisms
$S \mapsto L$
and $S \mapsto R$, respectively. The $\alpha _L$-classes of $S$ are
\[\{ (a_1, b_1); (a_1, b_2)\}\quad \hbox{and}\quad \{ (a_2, b_1); (a_2, b_2)\}.\] The
$\alpha _R$-classes of $S$ are
\[\{ (a_1, b_1); (a_2, b_1)\}\quad \hbox{and}\quad \{ (a_1, b_2); (a_2, b_2)\}.\]
It is easy to see that the congruences of $S$ are
$\iota_S$, $\alpha _L$, $\alpha _R$ and $\omega _S$.
We show that the ideals of ${\mathbb F}[S]$ are
\[{\mathbb F}[S], {\mathbb F}[\omega _S]=\{ \sum _{i,j=1}^2\alpha _{i,j}(a_i, b_j):\ \sum _{i,j=1}^2\alpha _{i,j}=0\},\]
\[J_L={\mathbb F}[\alpha _L], J_R={\mathbb F}[\alpha _R], J_L\cap J_R, \{ 0\}.\]
We note that \[dim({\mathbb F}[\omega _S])=3, dim(J_L)=dim(J_R)=2.\]

First we show that $J\subseteq {\mathbb F}[\omega _S]$ or $J={\mathbb F}[S]$ for every ideal $J$ of ${\mathbb F}[S]$.
Let $J$ be an arbitrary ideal of ${\mathbb F}[S]$. Assume that there is an
element
\[A=\alpha _{1,1}(a_1,b_1)+\alpha _{1,2}(a_1,b_2)+\alpha _{2,1}(a_2,b_1)+\alpha _{2,2}(a_2,b_2)\in J\] such that
$A\notin {\mathbb F}[\omega _S]$, that is
$\sum _{i,j=1}^2\alpha _{i,j}\neq 0$. Let $i, j\in \{ 1, 2\}$ be arbitrary
elements. Then
\[(\sum _{i,j=1}^2\alpha _{i,j})(a_i ,b_j)=(a_i, b_1)A(a_1, b_j)\in J.\] As $\sum _{i,j=1}^2\alpha _{i,j}\neq 0$, we get
$(a_i, b_j)\in J$ from which it follows that $S\subseteq J$. Consequently,
$J={\mathbb F}[S]$. Thus ${\mathbb F}[\omega _S]$ is the only maximal ideal of
${\mathbb F}[S]$.

Next we show that $J_L\cap J_R$ is the only ideal of ${\mathbb F}[S]$ whose dimension is $1$.
Let
\[A=\alpha _{1,1}(a_1,b_1)+\alpha _{1,2}(a_1,b_2)+\alpha _{2,1}(a_2,b_1)+\alpha _{2,2}(a_2,b_2)\in J_L\cap J_R\]
be an arbitrary element. As
\[(a_1, b_1)\ \alpha _L\ (a_1, b_2)\quad \hbox{and}\quad (a_2, b_1)\ \alpha _L\ (a_2, b_2),\] we have
\[\alpha _{1, 2}=-\alpha _{1, 1}\quad \hbox{and}\quad \alpha _{2, 2}=-\alpha _{2, 1}.\]
As
\[(a_1, b_1)\ \alpha _R\ (a_2, b_1)\quad \hbox{and}\quad (a_1, b_2)\ \alpha _R\ (a_2, b_2),\] we have
\[\alpha _{2, 1}=-\alpha _{1, 1}\quad \hbox{and}\quad \alpha _{2, 2}=-\alpha _{1, 2}.\]
Thus
\[A=\alpha _{1, 1} ((a_1,b_1)-(a_1,b_2)-(a_2,b_1)+(a_2,b_2))\]
for some $\alpha \in {\mathbb F}$. Consequently, the ideal $J_L\cap J_R$ is generated by
\[(a_1,b_1)-(a_1,b_2)-(a_2,b_1)+(a_2,b_2).\] Hence the dimension of $J_L\cap J_R$ is $1$.

To show that $J_L\cap J_R$ is the only ideal of ${\mathbb F}[S]$ whose dimension is $1$,
consider an ideal $J$ of ${\mathbb F}[S]$ generated by an element
\[0\neq B=\alpha _{1,1}(a_1,b_1)+\alpha _{1,2}(a_1,b_2)+\alpha _{2,1}(a_2,b_1)+\alpha _{2,2}(a_2,b_2).\]
Then $J\subset {\mathbb F}[\omega _S]$ and
\[(a_1, b_1)B=(\alpha _{1, 1}+\alpha _{2, 1})(a_1, b_1)+(\alpha _{1, 2}+\alpha _{2, 2})(a_1, b_2)\in J.\] Thus
there is a coefficient $\xi \in {\mathbb F}$ such that \[(a_1, b_1)B=\xi B.\]
Assume $\xi \neq 0$. Then $\alpha _{2, 1}=\alpha _{2, 2}=0$ and so
\[B=\alpha _{1, 1}(a_1, b_1)+\alpha _{1, 2}(a_1, b_2).\]
From
\[(a_2, b_2)B=\alpha _{1, 1}(a_2, b_1)+\alpha _{1, 2}(a_2, b_2)\in J\] we can
conclude that $\alpha _{1, 1}=\alpha _{1, 2}=0$ and so $B=0$. This is a contradiction.
Hence $\xi =0$.
Thus \[B=\alpha _{1, 1}(a_1, b_1)+\alpha _{1, 2}(a_1, b_2)\alpha _{1, 1}(a_2, b_1)-\alpha _{1, 2}(a_2, b_2).\]
As \[B(a_1, b_1)=(\alpha _{1,1}+\alpha _{1,2}((a_1, b_1)-(a_2, b_1))\in J,\]
we get $B(a_1, b_1)=\tau B$ for some $\tau \in {\mathbb F}$.
Assume $\tau \neq 0$. Then $\alpha _{1, 2}=0$ and so
\[B=\alpha _{1, 1}(a_1, b_1)-\alpha _{1,1}(a_2, b_1).\] From
\[B(a_2, b_2)=\alpha _{1,1}((a_1, b_2)-(a,2 ,b_2))\in J\] we can conclude that
$\alpha _{1, 1}=0$ and so $B=0$. This is a contradiction. Hence $\tau =0$.
Thus $\alpha _{1, 2}=-\alpha _{1, 1}$ and so
\[B=\alpha _{1, 1}((a_1, b_1)-(a_1, b_2)-(a_2, b_1)+(a_2, b_2))\in J_L\cap J_R.\]
As $J\neq \{ 0\}$ and $J_L\cap J_R$ is a minimal ideal of ${\mathbb F}[S]$, we get
\[J=J_L\cap J_R,\] that is,
$J_L\cap J_R$ is the only ideal of ${\mathbb F}[S]$ whose dimension is $1$.

As $dim(J_R+J_L)>dimJ_R$ and ${\mathbb F}[\omega _S]\supseteq J_R+J_L$, we have \[J_R+J_L={\mathbb F}[\omega _S].\]

Let $J$ be an arbitrary ideal of ${\mathbb F}[S]$ which differs from all of the ideals
${\mathbb F}[S], {\mathbb F}[\omega _S], J_L, J_R, J_L\cap J_R, \{ 0\}.$ Then $J\subset {\mathbb F}[\omega _S]$ and $dim(J)=2$.

If $J\cap J_L=\{ 0\}$ then $dim(J+J_L)=4$ which contradicts $J+J_L\subseteq {\mathbb F}[\omega _S]$.
Hence $dim(J\cap J_L)=1$ and so $J_L\cap J_R= J\cap J_L.$ From this we get $J\cap J_R=J_L\cap J_R.$
Recall that $C=(a_1, b_1)-(a_1, b_2)-(a_2, b_1)+(a_2, b_2)$ generates the ideal $J_L\cap J_R$. Let $A$ be an arbitrary element of $J\setminus (J_L\cap J_R)$. Then
$A$ and $C$ are linearly independent. So
\[A=\alpha(a_1, b_1)+\beta(a_1, b_2)+\gamma(a_2, b_1)+(-\alpha-\beta-\gamma)(a_2, b_2).\] If $\alpha=-\gamma$ then $A\in J_R$ which is a contradiction. Thus $\alpha\ne-\gamma$. Then
\[(a_1,b_1)A=(\alpha+\gamma)((a_1,b_1)-(a_1,b_2))\in J_L.\] As $J$ is an ideal and $A\in J$ we have
\[(a_1,b_1)A\in J\cap J_L=J_L\cap J_R.\] It means $\alpha=-\gamma$ which is also contradiction. Thus $Con(\mathbb{F}[S])$ is the next:

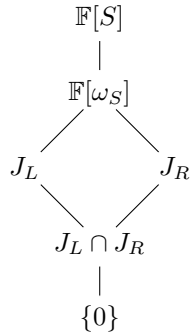
\begin{center}
\begin{tikzpicture}[scale=1,cap=round,
angle/.style={font=\fontsize{7}{7}\color{black}\ttfamily}
]
	
	\coordinate (a) at (0,-1);
	\coordinate (b) at (1,0);
	\coordinate (c) at (0,1);
	\coordinate (d) at (-1,0);
	\coordinate (e) at (0,2);
	\coordinate (f) at (0,-2);
	
	\draw (a) -- (b);
	\draw (a) -- (d);
	\draw (c) -- (b);
	\draw (d) -- (c);
	\draw (e) -- (c);
	\draw (f) -- (a);
			
	\filldraw[fill=white, draw=white]  (a) circle [radius=0.3];
	\filldraw[fill=white, draw=white]  (b) circle [radius=0.3];
	\filldraw[fill=white, draw=white]  (c) circle [radius=0.3];
	\filldraw[fill=white, draw=white]  (d) circle [radius=0.3];
	\filldraw[fill=white, draw=white]  (e) circle [radius=0.3];
	\filldraw[fill=white, draw=white]  (f) circle [radius=0.3];
		
	\draw (a) node{$J_L\cap J_R$};
	\draw (b) node{$J_R$};	
	\draw (c) node{$\mathbb{F}[\omega_S]$};
	\draw (d) node{$J_L$};
	\draw (e) node{$\mathbb{F}[S]$};
	\draw (f) node{$\{0\}$};

\end{tikzpicture}
 % \caption{}
\end{center}

It is a matter of checking to see that the $ker_{\varphi_{\{S;{\mathbb F}\}}}$-classes
of $Con({\mathbb F}[S])$ are
$\{\{0\},J_L\cap J_R\},\{J_L\},\{J_R\}$ and $\{\mathbb{F}[\omega_S],\mathbb{F}[S]\}$.
It is easy to see that $ker_{\varphi _{\{S;{\mathbb F}\}}}$ is $\vee$-compatible and so, by Lemma~\ref{xxx},
$\varphi _{\{S;{\mathbb F}\}}$ is a homomorphism of the semigroup
$(\mathrm{Con}({\mathbb F}[S]), \circ )$ onto the semigroup $(Con (S); \circ )$.
\hfill{\openbox}

\begin{corollary}\label{correctangularband} Let $S=L\times R$ be a rectangular band. Then, for a field $\mathbb F$,
$\varphi _{\{S;{\mathbb F}\}}$ is a homomorphism of the semigroup
$(\mathrm{Con}({\mathbb F}[S]), \circ )$ into the relation semigroup $({\cal B}_S;\circ)$
if and only if $|L|, |R|\leq 2$.
\end{corollary}

\noindent
{\bf Proof}. By Lemma~\ref{xxx}, Theorem~\ref{rectangularband} and Corollary 1.2 of \cite{Cherubini}, it is obvious.\hfill\openbox

\bigskip
%\vspace{-0.5cm}

\bigskip

\noindent
Department of Algebra

\noindent
Institute of Mathematics

\noindent
Budapest University of Technology and Economics

\noindent
1521 Budapest, Pf. 91

\noindent
e-mail (A. Nagy): nagyat@math.bme.hu

\noindent
e-mail (M. Zubor): zuborm@math.bme.hu

\end{document}